\documentclass[preprint,12pt]{elsarticle}

\usepackage{graphicx}
\usepackage{amssymb,amsmath}
\usepackage{tabularx,colonequals}
\usepackage{booktabs} 

\usepackage{amsthm}

\makeatletter
\let\saveqed\qed
\renewcommand\qed{%
  \ifmmode\displaymath@qed
  \else\saveqed
  \fi}
\makeatother

\usepackage{algorithm}
\usepackage{algorithmic}
\usepackage{lineno}

\newtheorem{theorem}{Theorem}
\newtheorem{assumption}{Assumption}
\newtheorem{proposition}[theorem]{Proposition}
\newtheorem{lemma}[theorem]{Lemma}

\newtheorem{example}{Example}
\newdefinition{definition}{Definition}
\newproof{pf}{Proof}

\def\eps{\varepsilon}
\def\tp{\tilde{p}}
\def\PKP{\text{PKP}}
\def\DP{\text{DP}}
\newcommand{\argmin}{\operatornamewithlimits{argmin}}
\newcommand{\argmax}{\operatornamewithlimits{argmax}}

\usepackage[usenames]{xcolor}
\newcommand{\edit}[1]{#1}
\begin{document}

\begin{frontmatter}


\title{The Product Knapsack Problem:\\ Approximation and Complexity}



\author[graz]{Ulrich Pferschy}
\author[graz]{Joachim Schauer}
\author[sr]{Clemens Thielen}

\address[graz]{Department of Statistics and Operations Research, University of Graz, Universitaetsstrasse~15, A-8010 Graz, Austria}

\address[sr]{TUM Campus Straubing, Technical University of Munich,\\ Am~Essigberg~3, D-94315~Straubing, Germany}

\begin{abstract}
We consider the product knapsack problem, which is the variant of the classical 0-1~knapsack problem where the objective consists of maximizing the product of the profits of the selected items.
These profits are allowed to be positive or negative.
We show that this recently introduced variant of the knapsack problem is weakly \textsf{NP}-hard and present a fully polynomial-time approximation scheme (FPTAS) for the problem. Moreover, we analyze the approximation quality achieved by a natural extension of the classical greedy procedure to the product knapsack problem.
\end{abstract}


\end{frontmatter}

% \linenumbers

\section{Introduction}\label{sec:introduction}

The 0-1~knapsack problem (KP) is a well-studied combinatorial optimization problem that has been treated extensively in the literature, with two monographs~\cite{MarTot90,Kellerer+etal:book} devoted to KP and its relatives. Given a positive knapsack capacity~$C$ and $n$~items $j=1,\dots,n$ with positive weights~$w_j$ and profits~$p_j$, the task in the classical 0-1 knapsack problem is to select a subset of items with maximum total profit subject to the constraint that the total weight of the selected items may not exceed the knapsack capacity. The 0-1 knapsack problem is $\textsf{NP}$-hard, but it admits a fully polynomial time approximation scheme (FPTAS) and can be solved exactly in pseudo-polynomial time by dynamic programming (cf.~\cite{Kellerer+etal:book}).

\newpage

The \emph{product knapsack problem (PKP)} is a new addition to the knapsack family. It has recently been introduced in~\cite{D-Ambrosio+etal:pkp} and is formally defined as follows:

\begin{definition}[Product Knapsack Problem (PKP)]\label{PKP-def}\mbox{}\\
  \begin{tabularx}{\linewidth}{lX}\hspace{-3mm}
    INSTANCE: & Items~$j\in N\colonequals \{1,\dots,n\}$ with weights~$w_j\in\mathbb{Z}$ and profits~$p_j\in\mathbb{Z}$, and a positive knapsack capacity~$C\in\mathbb{N}_+$. \\\hspace{-3mm}
    TASK: & Find a subset $S\subseteq N$ with $\sum_{j\in S}w_j\leq C$ such that $\prod_{j\in S} p_j$ is maximized.\
  \end{tabularx}
\end{definition}

The solution $S=\emptyset$ is always feasible and is assumed to yield an objective value of zero. \edit{Note that the assumption that the knapsack capacity as well as all weights and profits are integers is without loss of generality. Indeed, any instance with rational input data can be transformed into an equivalent instance with integer input data in polynomial time by multiplying all numbers by their lowest common denominator.}

\bigskip

D'Ambrosio et al.~\cite{D-Ambrosio+etal:pkp} list several application scenarios for PKP, in particular in the area of computational social choice, and also provide pointers to literature on other nonlinear knapsack problems. Furthermore, two different ILP formulations for PKP are presented and compared from both a theoretical and a computational perspective. In addition, D'Ambrosio et al.~\cite{D-Ambrosio+etal:pkp} develop an algorithm performing dynamic programming by weights with pseudopolynomial running time $\mathcal{O}(nC)$.
A computational study exhibits the strengths and weaknesses of the dynamic program and the ILP approaches for determining exact solutions of PKP depending on the characteristics of the test instances. However, no approximation results for PKP are known so far.

\subsection{Our Contribution}\label{subsec:contribution}

In this paper, we settle the complexity status of PKP by providing a proof of weak \textsf{NP}-hardness in Section~\ref{sec:complexity}, as well as an FPTAS based on dynamic programming by profits in Section~\ref{sec:fptas}.
The construction of an FPTAS deserves attention since standard greedy-type algorithms fail to give a constant approximation ratio for PKP. We demonstrate this in Section~\ref{sec:greedy} by providing a tight analysis of the greedy algorithm obtained by extending the classical greedy procedure for KP to PKP in the natural way.

\section{Preliminaries}\label{sec:prelim}

In contrast to KP, both the item weights~$w_j$ and the item profits~$p_j$ are allowed to be negative in PKP. However, one can exclude certain weight-profit combinations that yield ``useless'' items, which leads to the following assumption used throughout the paper:

\begin{assumption}\label{assum:useless-items}
Any instance of PKP satisfies:
\begin{enumerate}[(a)]
    \item Any single item fits into the knapsack, i.e., $w_j\leq C$ for all $j\in N$.
    \item All profits are nonzero, i.e., $p_j\in\mathbb{Z}\setminus \{0\}$ for all $j\in N$.
    \item For each item~$j\in N$ with negative profit~$p_j<0$, there exists another item~$j'\in N\setminus \{j\}$ with negative profit~$p_{j'}<0$ such that $w_j+w_{j'} \leq C$.
    \item All weights are nonnegative, i.e., $w_j\in\mathbb{N}_0$ for all $j\in N$.
    \item All items with weight zero have negative profit, i.e., $p_j<0$ whenever $w_j=0$.
\end{enumerate}
\end{assumption}

We note that Assumption~\ref{assum:useless-items} imposes no loss of generality and can easily be checked in polynomial time. Indeed, items~$j\in N$ violating~(a), (b), or~(c) can never be part of any feasible solution with positive objective value and may, thus, be removed from the instance. The nonnegativity of the weights~$w_j$ demanded in~(d) has been shown to impose no loss of generality in~\cite{D-Ambrosio+etal:pkp}. For~(e), we note that items~$j$ with $w_j=0$ can always be assumed to be packed if their profit is positive (but items~$j$ with $w_j=0$ and negative profit remain part of the optimization).

\bigskip

Using Assumption~\ref{assum:useless-items}~(b), the item set~$N$ can be partitioned into $N^+ \colonequals \{j \in N \mid p_j \geq 1\}$ and $N^- \colonequals \{j \in N \mid p_j \leq -1\}$. For convenience, we define
$p_{\max}\colonequals \max_{j\in N} |p_j|$, $p^+_{\max}\colonequals \max_{j\in N^+} p_j$, and $p^-_{\max}\colonequals \max_{j\in N^-} |p_j|$.

\bigskip

\edit{Throughout the paper,} we denote \edit{an optimal solution set for a given instance of PKP by~$S^*$} and the optimal objective value by~$z^*$. Note that we must always have $z^*\geq 1$ since packing any item from~$N^+$ or any feasible pair of items from~$N^-$ yields an objective value of at least~$1$.

\begin{definition}
For $0<\alpha\leq 1$, an algorithm~$A$ that computes a feasible solution set~$S\subseteq N$ with $\prod_{j\in S}p_j \geq \alpha\cdot z^*$ in polynomial time for every instance of PKP is called an \emph{$\alpha$-approximation algorithm} for PKP. The value~$\alpha$ is then called the \emph{approximation ratio} of~$A$.

\smallskip

A \emph{polynomial-time approximation scheme (PTAS)} for PKP is a family of algorithms~$(A_{\eps})_{\eps>0}$ such that, for each $\eps>0$, the algorithm~$A_{\eps}$ is a $(1-\eps)$-approximation algorithm for PKP. A PTAS~$(A_{\eps})_{\eps>0}$ for PKP is called a \emph{fully polynomial-time approximation scheme (FPTAS)} if the running time of~$A_{\eps}$ is additionally polynomial in $\frac{1}{\eps}$.
\end{definition}

Throughout the paper, $\log(x)$ always refers to the base~$2$ logarithm of~$x$ and $\ln(x)$ refers to the natural logarithm of~$x$.

\section{Complexity}\label{sec:complexity}

In this section, we show that PKP is weakly \textsf{NP}-hard.
In fact, \cite{D-Ambrosio+etal:pkp} already stated a short proof of weak \textsf{NP}-hardness as a side remark (Proposition~6).
Their proof, however, uses a reduction from KP and requires an exponential blow-up of the profits of the given instance of KP (by putting them into the exponent of 2). Since KP is only weakly \textsf{NP}-hard, this is problematic and does not prove the desired hardness result.

After publication of the first version of this manuscript in~\cite{PST19v1} in January 2019, we became aware that a proof of \textsf{NP}-hardness of PKP had already been published in~\cite{Halman:bicriteria} in December 2018. 
Since our proof uses only elementary operations and does not require the concepts of advanced calculus as in~\cite{Halman:bicriteria}, we present it here as a possibly useful alternative.

\bigskip

In our proof, we use a reduction from the \emph{Product Partition Problem (PPP)}, which was shown to be strongly \textsf{NP}-hard in~\cite{Ng+etal:product-partition} and can be stated as follows:

\begin{definition}[Product Partition Problem (PPP)]\label{PPP-def}\mbox{}\\
  \begin{tabularx}{\linewidth}{lX}\hspace{-3mm}
    INSTANCE: & Positive integers $a_1,\dots,a_n\in\mathbb{N}_+$.\\\hspace{-3mm}
    QUESTION: & Is there a subset $S\subseteq N=\{1,\dots,n\}$ such that\\ & $\prod_{j\in S} a_j = \prod_{j\in N\setminus S} a_j$?\
  \end{tabularx}
\end{definition}

\noindent
Within the proof, we use the following result:

\begin{lemma}\label{lem:log}
For all $x\geq 1$, we have
\begin{align*}
    \log(x+1) - \log(x) \geq \frac{1}{x+1}. 
\end{align*}
\end{lemma}

\begin{pf}
Considering the first derivative $(\log(x))'= \frac{1}{x \ln 2}$
and the concavity of the logarithm, we obtain
\begin{align*}
    \log(x) \leq \log(x+1) - \frac{1}{\ln 2 \, (x+1)} \leq \log(x+1) - \frac{1}{x+1}\,.\tag*{\qed}
\end{align*}
\end{pf}

\noindent
We are now ready to show our hardness result for PKP:

\begin{theorem}\label{th:nphard}
PKP is weakly \textsf{NP}-hard - even when all profits are positive.
\end{theorem}

\begin{pf}
Given an instance $a_1,\dots,a_n$ of PPP, we define 
\begin{align*}
M\colonequals (n+2)\cdot \left(\sqrt{\prod_{j=1}^n a_j}+1\right).
\end{align*}
Note that, since $\sqrt{\prod_{j=1}^n a_j}=\prod_{j\in S}a_j\in\mathbb{N}_+$ must hold for any solution~$S$ of the PPP instance, we can assume without loss of generality that $\sqrt{\prod_{j=1}^n a_j}\in\mathbb{N}$, which also implies that~$M$ is an integer.\footnote{The condition $\sqrt{\prod_{j=1}^n a_j}\in\mathbb{N}$ can be easily checked in polynomial time by using binary search for the square root within the interval $[1,\prod_{j=1}^n a_j]$.} Moreover, since we can put an index~$j$ with $a_j=1$ on either side of a partition $(S,N\setminus S)$ without changing the feasibility of~$S$ for the PPP instance, we can also assume that $a_j\geq 2$ for all $j\in N$.

\bigskip

We now construct an instance of PKP with positive profits as follows: There is one item for each index in $N=\{1,\dots,n\}$. Item~$j\in N$ has profit $p_j\colonequals a_j$ and weight $w_j\colonequals \lfloor M\cdot\log(a_j)\rfloor$. The knapsack capacity is set to $C\colonequals \lceil \frac{M}{2}\cdot \sum_{j=1}^n \log(a_j)\rceil$. %$C\colonequals \lceil M\cdot C'\rceil$, where $C'\colonequals \frac{1}{2}\cdot \sum_{j=1}^n \log(a_j)$.

\bigskip

Regarding the polynomial-time computability of the constructed PKP instance, we observe that the profits in the PKP instance clearly have polynomial encoding length in the size of the PPP instance. For the weights~$w_j$ and the knapsack capacity~$C$, we note that, while their definition involves logarithms~$\log(a_j)$ (which cannot be computed exactly in polynomial time), these logarithms only appear in expressions that are rounded to integers in the definition of the~$w_j$ and~$C$. Thus, the logarithms do not need to be computed exactly when computing the~$w_j$ and~$C$. The encoding lengths of the weights~$w_j$ and the knapsack capacity~$C$ are mainly determined by the encoding length of~$M$. Encoding~$M$ in binary requires $\Theta(\log(n)+\sum_{j=1}^n \log(a_j))$ bits. Thus, the binary encoded version of the constructed PKP instance has size polynomial in the size of the (binary or unary encoded) instance of PPP and can be constructed in polynomial time.\footnote{In contrast, encoding~$M$ in \emph{unary} requires $\Theta(n\sqrt{\prod_{j=1}^n a_j})$ bits, so the \emph{unary encoded} version of the PKP instance has size exponential in the size of both the binary and the unary encoded instance of PPP. Consequently, the reduction will indeed only show \emph{weak} \textsf{NP}-hardness, even though PPP is strongly \textsf{NP}-hard.}

\bigskip

We now show that the constructed instance of PKP has a solution with objective value at least $\sqrt{\prod_{j=1}^n a_j}$ if and only if the given instance of PPP is a yes-instance.

First assume that the given instance of PPP is a yes-instance, i.e., there exists $S\subseteq N$ such that $\prod_{j\in S} a_j = \prod_{j\in N\setminus S} a_j$. Then, as seen above, $\prod_{j\in S}p_j=\prod_{j\in S}a_j=\sqrt{\prod_{j=1}^n a_j}$, so packing exactly the items in~$S$ into the knapsack yields the desired profit. Moreover, this is a feasible solution for the PKP instance since raising both sides in $\prod_{j\in S}a_j=\sqrt{\prod_{j=1}^n a_j}$ to the $M$-th power and taking the logarithm yields that
\begin{align*}
   & \log\left(\left(\prod_{j\in S}a_j\right)^M\right)=\log\left(\left(\sqrt{\prod_{j=1}^n a_j}\right)^M\right) \\
    \Leftrightarrow\; & M\cdot \sum_{j\in S}\log(a_j)= \frac{M}{2} \cdot \sum_{j=1}^n\log(a_j) \\
    % \Leftrightarrow\; & \sum_{j\in S} M\cdot\log(a_j) =M\cdot C' \\
     % \Rightarrow\; & \sum_{j\in S} \lfloor M\cdot\log(a_j)\rfloor \leq \lceil M\cdot C'\rceil = C.
     \Rightarrow\; & \sum_{j\in S} \lfloor M\cdot\log(a_j)\rfloor \leq \left\lceil \frac{M}{2} \cdot \sum_{j=1}^n\log(a_j)\right\rceil = C.
\end{align*}

Conversely, assume that the constructed instance of PKP has a solution~$S$ with objective value at least $\sqrt{\prod_{j=1}^n a_j}$, i.e., $\prod_{j\in S}p_j = \prod_{j\in S}a_j \geq \sqrt{\prod_{j=1}^n a_j}$. We claim that we must then actually have equality, i.e., $\prod_{j\in S}a_j = \sqrt{\prod_{j=1}^n a_j}$, which directly implies that $\prod_{j\in S}a_j =\prod_{j\in N\setminus S}a_j$, so $S\subseteq N$ is a solution of the PPP instance.

\bigskip

Assume for the sake of a contradiction that $\prod_{j\in S}a_j > \sqrt{\prod_{j=1}^n a_j}$. Since both sides of the inequality are integers, this then implies that $\prod_{j\in S}a_j \geq \sqrt{\prod_{j=1}^n a_j}+1$, and raising both sides to the $M$-th power and taking the logarithm yields that
\begin{align}\label{eq:M-choice1}
   & \sum_{j\in S}M\cdot\log(a_j) \geq M\cdot \log\left(\sqrt{\prod_{j=1}^n a_j}+1\right).
\end{align}
In order to lower bound the right hand side of~\eqref{eq:M-choice1}, we use that, by
Lemma~\ref{lem:log},
\begin{align*}
\log\left(\sqrt{\prod_{j=1}^n a_j}+1\right)-\log\left(\sqrt{\prod_{j=1}^n a_j}\right)\geq \frac{1}{\sqrt{\prod_{j=1}^n a_j}+1},
\end{align*}
which implies that
\begin{align*}
    M & = (n+2)\cdot \left(\sqrt{\prod_{j=1}^n a_j}+1\right) \geq \frac{n+2}{\log(\sqrt{\prod_{j=1}^n a_j}+1)-\log(\sqrt{\prod_{j=1}^n a_j})}.
\end{align*}
Reordering terms, this yields
\begin{align}\label{eq:M-choice2}
    M\cdot \log\left(\sqrt{\prod_{j=1}^n a_j}+1\right) \geq M\cdot \log\left(\sqrt{\prod_{j=1}^n a_j}\right) +(n+2).
\end{align}
Combining~\eqref{eq:M-choice1} and~\eqref{eq:M-choice2} now shows that
\begin{align*}
   & \sum_{j\in S}M\cdot\log(a_j) \geq M\cdot \log\left(\sqrt{\prod_{j=1}^n a_j}\right)+(n+2) \\
   \Leftrightarrow & \sum_{j\in S}M\cdot\log(a_j) \geq \frac{M}{2} \cdot \sum_{j=1}^n\log(a_j) + (n+2) \\
   \Leftrightarrow & \sum_{j\in S}M\cdot\log(a_j) - n \geq \frac{M}{2} \cdot \sum_{j=1}^n\log(a_j) + 2
\end{align*}
Since rounding down each of the at most~$n$ summands on the left decreases the sum by at most~$n$ and rounding up $\frac{M}{2} \cdot \sum_{j=1}^n\log(a_j)$ increases this value by at most one, this implies that 
\begin{align*}
    \sum_{j\in S} \lfloor M\cdot\log(a_j)\rfloor \geq \left\lceil \frac{M}{2} \cdot \sum_{j=1}^n\log(a_j)\right\rceil + 1 = C+1,
\end{align*}
which contradicts the feasibility of the set~$S$ for the PKP instance and completes the proof.\qed
\end{pf}

\section{A Fully Polynomial-Time Approximation Scheme}\label{sec:fptas}

We now derive a fully polynomial-time approximation scheme (FPTAS) for PKP \edit{based on dynamic programming.

The most common approach for the exact solution of knapsack-type problems in pseu\-do-po\-ly\-no\-mial time applies dynamic programming by weights. 
This means that, for every capacity value $d=0,1,\ldots,C$, the largest profit value reachable by a feasible solution is determined, which yields a running time polynomial in~$C$ (see \cite[Sec.~2.3]{Kellerer+etal:book}).
However, for obtaining fully polynomial-time approximation schemes, one usually performs dynamic programming by profits.
In this case, for every profit value~$p$ up to some upper bound~$U$ on the objective function value, the smallest weight required for a feasible solution with profit~$p$ is sought, which leads to a running time polynomial in~$U$ (see \cite[Lemma~2.3.2]{Kellerer+etal:book}). Then, the profit space is simplified in some way, e.g., by scaling (cf.~\cite[Sec.~2.6]{Kellerer+etal:book}), such that the running time of the dynamic program becomes polynomial and the incurred loss of accuracy remains bounded.
} %end of edit

D'Ambrosio et al.~\cite{D-Ambrosio+etal:pkp} provide an algorithm \edit{solving PKP with} dynamic programming by weights, where each entry of the dynamic programming array contains the objective value of a subproblem.
However, exchanging the roles of profits and weights (as it is done, e.g., for KP, see \cite[Sec.~2.3]{Kellerer+etal:book}),
would require a dynamic programming array of length $\mathcal{O}(p_{\max}^n)$, which is exponential and does not permit a suitable scaling procedure.

An obvious way out of this dilemma would be the application of the logarithm to the profits. In fact, such an approach is suggested as a side remark in \cite[Sec.~3]{D-Ambrosio+etal:pkp} for dynamic programming by weights (without commenting on the details of the rounding process).
For dynamic programming by profits, however, the profit values must be mapped to integers as indices of the dynamic programming array and there seems to be no way to preserve optimality in such a process.
It should also be noted that applying any $k$-approximation algorithm for KP to the instance resulting from logarithmization would only yield a $(1/z^*)^{1/k}$-approximation for PKP.
Thus, constant-factor approximations for PKP require different approaches.

\bigskip

We now construct a scaled profit space that actually yields a $(1-\eps)$--approximation for PKP. Our scaling construction is based on a parameter $K>0$ depending on~$\eps$, which will be defined later. For every item~$j$, we define an integer scaled profit value in the logarithmized space as
\begin{align}
    \tp_j\colonequals \left\lfloor\frac{\log(|p_j|)}{K}\right\rfloor. \label{eq:scaledprofit}
\end{align}
Since $|p_j|\geq 1$, we have $\tp_j\geq 0$, and we obtain $\tp_j=0$ if and only if $|p_j|=1$.
Note that an item~$j$ with $p_j=-1$ and $\tp_j=0$ might still be useful for changing the sign of the solution of PKP.
Analogous to $p_{\max}$, we define 
$\tp_{\max} \colonequals \left\lfloor\frac{\log(p_{\max})}{K}\right\rfloor$.
Ruling out trivial instances, we can assume without loss of generality that $p_{\max}\geq 2$, so $\log(p_{\max})\geq 1$.
As discussed in the proof of Theorem~\ref{th:nphard}, also for the FPTAS we require the logarithm only for expressions which are rounded down to integers and, thus, do not have to compute these values exactly.

\bigskip

We define the following dynamic programming arrays for profit values $\tp=0,1,\ldots, n\cdot \tp_{\max}$:
\begin{eqnarray*}
W_j^+(\tp) &\colonequals& \min_{S\subseteq \{1,\ldots,j\}} \left\{\sum_{i\in S} w_i \mid \sum_{i \in S} \tp_i = \tp,\: |S\cap N^-| \mbox{ is even} \right\},\\
W_j^-(\tp) &\colonequals& \min_{S\subseteq \{1,\ldots,j\}} \left\{\sum_{i\in S} w_i \mid \sum_{i \in S} \tp_i = \tp,\: |S\cap N^-| \mbox{ is odd} \right\}.
\end{eqnarray*}
Note that the empty set has even cardinality. For convenience, we set the minimum over the empty set equal to~$+\infty$.

\bigskip

The computation of these arrays can be done by the following recursion, which is related to Algorithm $\DP_{\PKP}$ in~\cite[Fig.~1]{D-Ambrosio+etal:pkp}:
\begin{align*}%\arrayrowsep=1.4pt\def\arraystretch{2.2}
\begin{array}{ll}
    \mbox{If } p_j\geq 1, \mbox{ then: } \\
    & W_j^+(\tp)\colonequals \min \{W_{j-1}^+(\tp),\, W_{j-1}^+(\tp-\tp_j) +w_j\} \\ \addlinespace
    & W_j^-(\tp)\colonequals \min \{W_{j-1}^-(\tp),\, W_{j-1}^-(\tp-\tp_j) +w_j\} \\
    \mbox{If } p_j\leq -1, \mbox{ then: } \\
    & W_j^+(\tp)\colonequals \min \{W_{j-1}^+(\tp),\, W_{j-1}^-(\tp-\tp_j) +w_j\} \\ \addlinespace
    & W_j^-(\tp)\colonequals \min \{W_{j-1}^-(\tp),\, W_{j-1}^+(\tp-\tp_j) +w_j\} 
\end{array}
\end{align*}
The obvious initialization is given by $W_0^+(0)\colonequals0$
and setting all other entries (including the hypothetical ones with $\tp<0$) to~$+\infty$.

\bigskip

The approximate solution set~$S^A$ is represented by the array entry with 
$\max \{\tp \mid W_n^+(\tp) \leq C\}$.
It follows by construction that $S^A$ maximizes the total profit in the associated instance of KP with scaled profits~$\tp_j$ among all subsets of~$N$ that fulfill the weight restriction and contain an even number of items from~$N^-$. 

\bigskip

\noindent
In the following, we show that, by choosing
\begin{align}\label{eq:defis}
  K\colonequals\frac{\eps}{n^2}>0\,,
\end{align}
the set~$S^A$ yields a $(1-\eps)$-approximation for PKP and can be computed in polynomial time via the above dynamic programming procedure. To this end, we use the following two lemmas: 

\begin{lemma}\label{lem:epsi} For $\eps\in(0,1)$, we have $\eps \leq - \log (1-\eps)$.
\end{lemma}

\begin{pf}
The statement follows since, for any $x\in(0,1)$, we have
\begin{align*}
    -\log(1-x) = -\ln(1-x) / \ln 2 \geq -\ln(1-x) = \sum_{k=1}^\infty \frac{x^k}{k} \geq x.\tag*{\qed}
\end{align*}
\end{pf}

\begin{lemma}\label{lem:at-least-p-max}
Any optimal solution set~$S^*$ for PKP satisfies
\begin{align*}
    \sum_{j\in S^*}\log(|p_j|) \geq \log(p_{\max}).
\end{align*}
\end{lemma}

\begin{pf}
Let~$j_{\max}\in N$ denote an item with $|p_{j_{\max}}|=p_{\max}$.
If $p_{j_{\max}}>0$, then the set $\{j_{\max}\}$, which is feasible for PKP by Assumption~\ref{assum:useless-items}~(a), has objective value~$p_{\max}$. If $p_{j_{\max}}<0$, Assumption~\ref{assum:useless-items}~(c) implies that there exists another item $j'\neq j_{\max}$ with $p_{j'}<0$ such that $\{j_{\max},j'\}$ is feasible for PKP, and this set has objective value $p_{j'}\cdot p_{j_{\max}} \geq p_{\max}$ since $p_{j'}\leq-1$ by Assumption~\ref{assum:useless-items}~(b). Thus, in both cases, the optimality of~$S^*$ for PKP implies that
\begin{align*}
   & \prod_{j\in S^*}p_j \geq p_{\max} \\
   \Leftrightarrow\; & \log\left(\prod_{j\in S^*}p_j\right) \geq \log(p_{\max}) \\
   \Leftrightarrow\; & \log\left(\prod_{j\in S^*}|p_j|\right) \geq \log(p_{\max}) \\
   \Leftrightarrow\; & \sum_{j\in S^*}\log(|p_j|) \geq \log(p_{\max}).\tag*{\qed}
\end{align*}
\end{pf}

\begin{proposition}\label{prop:fptastime}
The running time for computing~$S^A$ is in 
$\mathcal{O}(\frac{n^4}{\eps} \log(p_{\max}))$,
which is polynomial in $1/\eps$ and the encoding length of the input of PKP.
\end{proposition}
\begin{pf}
Clearly, for each of the $n$~items, one has to pass through the whole length of the two dynamic programming arrays.
Therefore, the total running time is in
\begin{align*}
    \mathcal{O}(n^2\,\tp_{\max}) = \mathcal{O}\left(n^2 \frac{\log(p_{\max})}{K}\right) = \mathcal{O}\left(n^4 \frac{\log(p_{\max})}{\eps}\right).\tag*{\qed}
\end{align*}
\end{pf}

%\newpage

\begin{proposition}\label{prop:fptaserror}
The set~$S^A$ yields a $(1-\eps)$--approximation for PKP.
\end{proposition}

\begin{pf}
%Let~$S^*$ denote an optimal solution set of PKP.

%\bigskip

The proof consists of two parts. First, we analyze the effect of scaling by~$K$ and rounding down in \eqref{eq:scaledprofit} by showing that~$S^A$ yields an objective value close to the value of \edit{an optimal solution set}~$S^*$ for \edit{PKP in} the associated instance of KP with profits $\log(|p_j|)$. The argumentation closely follows the standard FPTAS for KP (see \cite[Sec.~2.6]{Kellerer+etal:book}):
\begin{align}
\sum_{j\in S^A} \log(|p_j|) &\geq
\sum_{j\in S^A} K \cdot \left\lfloor\frac{\log(|p_j|)}{K}\right\rfloor\label{eq:chain-first}\\
&\geq  
\sum_{j\in S^*} K \cdot \left\lfloor\frac{\log(|p_j|)}{K}\right\rfloor\label{eq:optsol}\\
&\geq  
\sum_{j\in S^*} K \cdot \left(\frac{\log(|p_j|)}{K}-1\right)\\
&\geq
\sum_{j\in S^*} \log(|p_j|) - n\cdot K \label{eq:chain-last}
\end{align}
In~\eqref{eq:optsol}, we exploited the optimality of~$S^A$ for the KP instance with profits~$\tp_j$. We now set
\begin{align}\label{eq:eps-def}
  \eps' \colonequals \frac{-\log(1-\eps)}{n\cdot\log(p_{\max})}>0.
\end{align}
Then, using the definition of~$K$ in~\eqref{eq:defis} and that $\eps \leq - \log (1-\eps)$ for $\eps\in(0,1)$, we obtain
\begin{align*}
    n\cdot K = \frac{\eps}{n} \leq \frac{-\log(1-\eps)}{n} = \eps'\cdot \log(p_{\max}),
\end{align*}
and using that $\sum_{j\in S^*}\log(|p_j|) \geq \log(p_{\max})$ by Lemma~\ref{lem:at-least-p-max}, the chain of inequalities in \eqref{eq:chain-first}--\eqref{eq:chain-last} yields that 
\begin{align*}
\sum_{j\in S^A} \log(|p_j|) \geq 
\sum_{j\in S^*} \log(|p_j|) - \eps'\cdot \log(p_{\max})
\geq (1-\eps') \sum_{j\in S^*} \log(|p_j|).
\end{align*}
In the second part of the proof, we simply raise two to the power of both sides of this inequality, i.e., $2^{\sum_{j\in S^A} \log(|p_j|)} \geq 
\left(2^{\left(\sum_{j\in S^*} \log(|p_j|)\right)}\right)^{1-\eps'}$, so
\begin{align}
\prod_{j\in S^A} |p_j| &\geq 
\left(\prod_{j\in S^*} |p_j|\right)^{1-\eps'} \\
&= z^* \cdot (1/z^*)^{\eps'}\\
\; &\geq z^* \cdot \left(\frac{1}{(p_{\max})^n}\right)^{\eps'}\label{eq:optupper}\\
\; &= z^* \cdot 2^{-\eps'\, n \log(p_{\max})}\\
\; &= z^* \cdot 2^{\log(1-\eps)} = (1-\eps) z^* \label{eq:eps-plugin}
\end{align}
Here, \eqref{eq:optupper} \edit{is derived} from the trivial upper bound~$z^*\leq (p_{\max})^n$, and~\eqref{eq:eps-plugin} from the definition of~$\eps'$ in~\eqref{eq:eps-def}.
Recalling that~$S^A$ contains an even number of items from $N^-$, the statement follows.\qed
\end{pf}

\noindent
Propositions~\ref{prop:fptastime} and~\ref{prop:fptaserror} immediately yield the following theorem:

\begin{theorem}\label{thm:fptas}
There exists an FPTAS with running time in $\mathcal{O}(\frac{n^4}{\eps}  \log(p_{\max}))$ for PKP.\qed
\end{theorem}

\edit{
\noindent
We conclude this section with an example illustrating how the FPTAS works.

\begin{example}\label{ex:fptas}
Consider the instance of PKP given by the $n=5$~items with profits and weights as shown in Table~\ref{tab:examplefptasprof} and a knapsack capacity of $C \colonequals 9$. We choose~$\eps = 0.025$ so that $K=\frac{\eps}{n^2}=0.001$.

\begin{table}[ht]
\edit{
\centering
\setlength{\tabcolsep}{2mm}
\begin{tabular}{|r|ccccc|}\hline
	item~$j$ & 1 & 2 & 3 & 4 & 5 \\\hline
	$p_j$ & $1$ & $2^{10}-1$ & $-(2^{10}+1)$ & $2^{10}$ & $-1$\\
	$w_j$ & $1$ & $5$ & $5$ & $5$ & $4$ \\
	$\tilde{p}_j$ & $0$ & $9998$ & $10001$ & $10000$ & $0$\\
	\hline
\end{tabular}\label{tab:examplefptasprof}
\caption{Profits~$p_j$, weights~$w_j$, and scaled profits~$\tp_j$ of the items in Example~\ref{ex:fptas}.} 
} % end of \edit
\end{table}

The resulting scaled profits~$\tp_j$ are shown in the last row of Table~\ref{tab:examplefptasprof} and we have $\tp_{\max}=10001$, so $n\cdot\tp_{\max}=50005$.
Thus, the FPTAS computes the relevant dynamic programming arrays~$W_j^+(\tp)$ and~$W_j^-(\tp)$ for all profit values $\tp=0,1,\ldots, 50005$. Note that $N^+=\{1,2,5\}$ and~$N^-=\{3,5\}$.

For this instance, the FPTAS finds the optimal solution set~$S^*=\{3,5\}$ during the computation of~$W_5^+(10001)$, which is given as follows: 
\begin{align*}
    W_5^+(10001) = \min \left\{W_{4}^+(10001),\, W_{4}^-(10001-0) + 4\right\} = \min\left\{+\infty, 9\right\}=9
\end{align*}

Here, $W_{4}^+(10001) = +\infty$ since a scaled profit of~$10001$ cannot be obtained by any subset~$S\subseteq\{1,2,3,4\}$ containing an even number of items from~$N^-$, and $W_{4}^-(10001) = 5$ since a scaled profit of~$10001$ is reachable by the subset~$S=\{3\}\subseteq\{1,2,3,4\}$ that contains an odd number of items from~$N^-$. Thus, the solution set corresponding to the array entry~$W_5^+(10001)$ is $\{3,5\}=S^*$, and since~$\tp=10001$ is indeed the highest value of~$\tp$ for which $W_5^+(\tp)\leq C=9$, this is also the set~$S^A$ returned by the FPTAS.
\end{example}
}

\section{A Greedy-Type Algorithm}\label{sec:greedy}
For KP, the classical greedy procedure is probably one of the most obvious first attempts for anybody confronted with the problem. Hence, it is interesting to evaluate the performance of a variant of this greedy procedure also for PKP.

It is known that, for obtaining a bounded approximation ratio for KP in the classical greedy procedure, one has to take into account also the item with largest profit as a singleton solution (cf.~\cite[Sec.~2.1]{Kellerer+etal:book}). Extending this requirement to the negative profits allowed in PKP, we additionally determine, among all items with negative profits, a feasible pair of items with largest profit product. Moreover, if the greedy solution contains an odd number of items from $N^-$, we simply remove the negative-profit item \edit{whose profit has the smallest absolute value}. This leads to the following natural greedy algorithm for PKP, which we refer to as \textsc{Product Greedy}:

\begin{algorithm}[H]
\begin{algorithmic}[1]
\STATE Sort and renumber the items in nonincreasing order of $\frac{\log(|p_j|)}{w_j}$.\\ (Items~$j$ with $w_j=0$ are put to the beginning of the ordering.)
 \label{alg:1} 
\STATE Perform the classical greedy procedure with this ordering yielding solution set $\edit{\bar{S}} \subseteq N$.
\IF{$|\edit{\bar{S}} \cap N^-|$ is odd}   
\STATE $j^- \colonequals \argmin \{|p_j| \mid j \in \edit{\bar{S}} \cap N^-\}$
\STATE $S \colonequals \edit{\bar{S}}\setminus \{j^-\}$
\ELSE
\STATE \edit{$S \colonequals \bar{S}$}
\ENDIF
\STATE Let $\{j_1, j_2\} \subseteq N^-$ be a pair of items with $w_{j_1}+w_{j_2} \leq C$ maximizing the profit product $p_{j_1}\cdot p_{j_2}$ over all such pairs. 
%\STATE Let $\{j_1, j_2\}$ be a pair of items with maximal multiplied profits,
%whose profit values have the same sign and $w_{j_1}+w_{j_2} \leq C$. 
\STATE Let $j^+_{\max} \colonequals \argmax \{p_j \mid j \in N^+\}$ be an item with largest positive profit.
\RETURN the best among the three solutions $S$, $\{j_1,j_2\}$, and $\{j^+_{\max}\}$.
% \STATE $z^H \colonequals \max\{\prod_{j\in S} p_j,\, p_{j_1} \cdot p_{j_2},\, p^+_{\max} \}$
\end{algorithmic}
\caption{Algorithm \textsc{Product Greedy}}
\end{algorithm}

We note that, since $\log(|p_j|)/w_j=\log\left(|p_j|^{1/w_j}\right)$ and the logarithm is a strictly increasing function, the sorting and renumbering of the items in step~$1$ of \textsc{Product Greedy} can equivalently be done by sorting the items in nonincreasing order of $|p_j|^{1/w_j}$, which means that the values $\log(|p_j|)/w_j$ do not have to be computed in the algorithm.

We let $j^+_{\max} \colonequals \argmax \{p_j \mid j \in N^+\}$ denote an item with largest positive profit (i.e., with $p_{j^+_{\max}} = p^+_{\max}$) as in \textsc{Product Greedy}. Similarly, we let $j^-_{\max} \colonequals \argmax \{|p_j| \mid j \in N^-\}$ denote an item with smallest negative profit \edit{(i.e., with $-p_{j^-_{\max}}=p^-_{\max}$)}. Then, by Assumption~\ref{assum:useless-items}~(c), there exists another item in~$N^-$ that can be packed into the knapsack together with~$j^-_{\max}$. This implies that the profits of the items~$j^-$ and $j_1,j_2$ considered in \textsc{Product Greedy} satisfy
\begin{align}
    p_{j_1} \cdot p_{j_2} \geq -p_{j^-_{\max}} \geq -p_{j^-}. \label{eq:negpair}
\end{align}

In the following analysis, we denote the objective value obtained by \textsc{Product Greedy} by~$z^H$.

\begin{theorem}\label{th:greedy}\mbox{}\\
\begin{enumerate}[(a)]\vspace{-7mm}
    \item \textsc{Product Greedy} is a $1/(z^*)^{2/3}$-approximation algorithm for PKP.
    \item \textsc{Product Greedy} is a $1/(p_{\max})^2$-approximation algorithm for PKP.
\end{enumerate}
\end{theorem}

\begin{pf}
The algorithm clearly runs in polynomial time. In order to analyze its approximation ratio, let $s\in N$ be the \emph{split item}, i.e., the first item in the given order that cannot be packed into the knapsack anymore during the greedy procedure performed in step~$2$. Similar to the analysis of the greedy procedure for KP, the analysis concentrates on the \emph{split solution}, i.e., the set of items \edit{$\bar{S}= \{j \in N \mid j \leq s-1\}$ produced in step~$2$ of \textsc{Product Greedy}}.

% Again let~$S^*$ denote an optimal solution set of PKP.
We distinguish two cases depending on the number of items with negative profits in~$\edit{\bar{S}}$ and, for each of the two cases, two subcases depending on the sign of the profit~$p_s$ of the split item~$s$:

\medskip

\noindent
Case~1: $|\edit{\bar{S}} \cap N^-|$ is even.\\
\edit{In this case, the solution~$S=\bar{S}$ is considered when choosing the best solution in step~$11$.} Consider the sign of the split item's profit. If $p_s>0$, then 
\begin{align*}
2\cdot \log(z^H)
& \geq 2\cdot \max\left\{\sum_{j \in \bar{S}} \log(|p_j|), \, \log(p^+_{\max})\right\} \\
& \geq \sum_{j \in \bar{S}} \log(|p_j|) + \log(p^+_{\max}) \\
& \geq \sum_{j=1}^s \log(|p_j|).
\end{align*}
Obviously, we also have $\log(z^H) + \log(p^+_{\max})  \geq \sum_{j=1}^s \log(|p_j|)$.

\bigskip

\noindent
Similarly, if $p_s<0$, then 
\begin{align*}
2\cdot \log(z^H)
& \geq 2\cdot \max\left\{\sum_{j \in \bar{S}} \log(|p_j|), \, \log(|p_{j_1}|)+\log(|p_{j_2}|)\right\} \\
& \geq \sum_{j \in \bar{S}} \log(|p_j|) + \log(|p_{j_1}|)+\log(|p_{j_2}|) \\
& \geq \sum_{j \in \bar{S}} \log(|p_j|) + \log(|p_s|) \\
& = \sum_{j=1}^s \log(|p_j|),
\end{align*}
where the third inequality follows from~\eqref{eq:negpair}. Moreover, we have
$\log(z^H) + \log(p^-_{\max}) \geq \sum_{j=1}^s \log(|p_j|)$.

% \bigskip
%\newpage

\noindent
Case~2: $|\edit{\bar{S}} \cap N^-|$ is odd.\\
In this case, \edit{the solution $S=\bar{S}\setminus\{j^-\}$ is considered when choosing the best solution in step~$11$}. If $p_s>0$, we obtain
\begin{align*}
3\cdot \log(z^H)
& \geq 3 \cdot\max\left\{\sum_{j \in \bar{S}\edit{\setminus\{j^-\}}} \log(|p_j|), \, \log(|p_{j_1}|)+\log(|p_{j_2}|),\, \log(p^+_{\max})\right\} \\
& \geq \sum_{j \in \bar{S}\edit{\setminus\{j^-\}}} \log(|p_j|) + \log(|p_{j_1}|)+\log(|p_{j_2}|) + \log(p^+_{\max}) \\
& \geq \sum_{j \in \bar{S}\edit{\setminus\{j^-\}}} \log(|p_j|) + \log(|p_{j^-}|) + \log(p_s) \\
& \edit{=} \sum_{j=1}^s \log(|p_j|),
\end{align*}
\edit{by} invoking \eqref{eq:negpair} again. In this case, we also have
$\log(z^H) + \log(p^-_{\max}) + \log(p^+_{\max})  \geq \sum_{j=1}^s \log(|p_j|)$.

\bigskip

\noindent
Similarly, if $p_s<0$, then
\begin{align*}
3\cdot\log(z^H)
& \geq 3 \cdot\max\left\{\sum_{j \in \bar{S}\edit{\setminus\{j^-\}}} \log(|p_j|), \, \log(|p_{j_1}|)+\log(|p_{j_2}|)\right\} \\
& \geq \sum_{j \in \bar{S}\edit{\setminus\{j^-\}}} \log(|p_j|) + 2(\log(|p_{j_1}|)+\log(|p_{j_2}|))) \\
& \geq \sum_{j \in \bar{S}\edit{\setminus\{j^-\}}} \log(|p_j|) + \log(|p_{j^-}|) + \log(|p_s|) \\
& \edit{=} \sum_{j=1}^s \log(|p_j|).
%\geq \sum_{j\in S^*} \log(|p_j|).
\end{align*}
Moreover, we have $\log(z^H) + 2\log(p^-_{\max}) \geq \sum_{j=1}^s \log(|p_j|)$.

\bigskip

\noindent
Summarizing all four cases, we always have
\begin{align*}
    3\cdot \log(z^H) \geq \sum_{j=1}^s \log(|p_j|).
\end{align*}
%Now let~$S^*$ be an optimal solution set of PKP.
Then, since $\sum_{j=1}^s \log(|p_j|)$ is 
an upper bound on the optimal objective value of the LP relaxation of the associated instance of KP with profits $\log(|p_j|)$ (see, e.g.,~\cite{Kellerer+etal:book}), we have $\sum_{j=1}^s \log(|p_j|) \geq \sum_{j\in S^*} \log(|p_j|)
=\log(\prod_{j\in S^*} p_j)=\log(z^*)$ (clearly, $|S^* \cap N^-|$ must be even).
This yields
\begin{align*}
    3\cdot \log(z^H) \geq \log(z^*) \Longleftrightarrow z^H \geq (z^*)^{1/3}
\end{align*}
and proves the approximation ratio in~(a).

Moreover, in all four cases the additive error in the logarithmic space can be bounded by 
$\max\{\log(p^+_{\max}), \, \log(p^-_{\max})\} + \log(p^-_{\max}) \leq 2\cdot \edit{\log}(p_{\max})$, which yields the approximation ratio in~(b).\qed
\end{pf}

The approximation ratios obtained by \textsc{Product Greedy} are rather disappointing. The following example, however, shows that the analysis in the proof of Theorem~\ref{th:greedy} is asymptotically tight and that a considerable deviation from the greedy principle would be necessary to improve upon the obtained approximation ratios:

\begin{example}\label{ex:greedy-tight}
Consider the instance of PKP given by the item profits and weights shown in Table~\ref{tab:exampledata} and a knapsack capacity of $C\colonequals 3M$ for some large, positive integer~$M$.

\begin{table}[H]
\centering
\setlength{\tabcolsep}{2mm}
\begin{tabular}{|r|cccccc|}\hline
	item~$j$ & 1 & 2 & 3 & 4 & 5 & 6 \\\hline
	$p_j$ & $2$ & $M+2$ & $-(M+1)$ & $M$ & $M$ & $-1$\\
	$w_j$ & $1$ & $M$ & $M$ & $M$ & $M$ & $M$\\
	\hline
\end{tabular}\label{tab:exampledata}
\caption{Profits~$p_j$ and weights~$w_j$ of the items in Example~\ref{ex:greedy-tight} with items indexed in nonincreasing order of $\log(|p_j|)/w_j$.}
\end{table}

Algorithm \textsc{Product Greedy} first finds $\edit{\bar{S}}=\{1,2,3\}$ in step~$2$, but has to remove item~$3$ in step~$5$ since $|\edit{\bar{S}}\cap N^-|=1$, which yields $S=\{1,2\}$ with an objective value of~$2(M+2)$. The best negative pair found in step~$\edit{9}$ is given by~$j_1=3$ and~$j_2=6$, and has profit product~$M+1$. Finally, $j^+_{\max}=2$ with $p_{j^+_{\max}}=p^+_{\max}=M+2$ in step~$\edit{10}$. Therefore, \textsc{Product Greedy} returns the solution $\{1,2\}$ with an objective value of $z^H=2(M+2)$, while the optimal solution consists of items \edit{$2,4$, and~$5$} with objective value
\edit{$z^*=(M+2)M^2$}.
\end{example}



%% References


% \section*{References}
\bibliographystyle{elsarticle-num}
\bibliography{pkp.bib}






% \end{thebibliography}

\end{document}